\documentclass[reqno]{amsart}

\usepackage[applemac]{inputenc}

\usepackage{amsmath}
\usepackage{amssymb}
\usepackage{amsfonts}
\usepackage{graphicx}
\usepackage{amsthm}
\usepackage{enumerate}
\usepackage{lscape}
\usepackage{dsfont}
\usepackage{color}
\usepackage{mathtools}

\usepackage{setspace}
\newcommand{\R}{\mathds{R}}                   
\newcommand{\Z}{\mathds{Z}}

\newcommand{\CP}{\mathds{C}\mathrm{P}}

\newcommand{\Ric}{\mathrm{Ric}}
                       
\newcommand{\f}{\rightarrow}                  
\newcommand{\C}{\mathds{C}}            
\newcommand{\de}{\partial}          

\newcommand{\K}{K\"{a}hler}

\newcommand{\ov}[1]{\overline{#1}}

\newcommand{\deb}{\ov\partial}

\newcommand{\Cut}{{\operatorname{Cut}}}
\newcommand{\Vol}{{\operatorname{Vol}}}
\newcommand{\cat}{{\operatorname{cat}}}

\newcommand{\ep}{\epsilon} 

\newcommand{\rank}{\operatorname{rank}}

\newcommand{\I}{\operatorname{I}}
\newcommand{\II}{\operatorname{II}}
\newcommand{\III}{\operatorname{III}}
\newcommand{\IV}{\operatorname{IV}}

\newtheorem{theor}{Theorem}

\newtheorem{lem}[theor]{Lemma}
\newtheorem{cor}[theor]{Corollary}

\newtheorem{remark}[theor]{Remark}

\newtheorem*{theorA}{Theorem A}

\begin{document}

\title[Minimal symplectic atlases of Hermitian symmetric spaces]{Minimal symplectic atlases of   Hermitian symmetric spaces}

\author[R. Mossa]{Roberto Mossa}
\address[R. Mossa]{Dipartimento di Matematica e Informatica \\
         Universit\`a degli studi di Cagliari (Italy)}
         \email{roberto.mossa@gmail.com}

\author[G. Placini]{Giovanni Placini}
\address[G. Placini]{Dipartimento di Matematica \\
         Universit\`a di Pisa (Italy)}
         \email{giovanniplacini@tiscali.it}

\thanks{
The first  author was  supported by Prin 2010/11 -- Variet\`a reali e complesse: geometria, topologia e analisi armonica -- Italy}
\subjclass[2010]{53D05;  53C55;  53D05; 53D45} 
\keywords{Minimal symplectic atlases; Darboux chart; Gromov width; Hermitian symmetric spaces of compact type.}

\begin{abstract}
In this paper we estimate the minimal number of Darboux charts needed to cover a Hermitian symmetric space of compact type $M$ in terms of the degree of their embeddings in $\C P^N$. The proof is based on the recent work of Y. B. Rudyak and F. Schlenk \cite{minatlas} and on the symplectic geometry tool developed by the first author in collaboration with A. Loi and F. Zuddas \cite{GWHSSCT}. As application we compute this number for a large class of  Hermitian symmetric spaces of compact type.
\end{abstract}
 
\maketitle

\section{Introduction and statements of the main results}
Consider the open ball of radius $r$,
\begin{equation*}\label{ball}
B^{2n}(r)=\{(x, y)\in\R^{2n}\  |\  \sum_{j=1}^nx_j^2+y_j^2<r^2 \}
\end{equation*}
in the standard symplectic space $(\R^{2n}, \omega_0)$, where $\omega_0=\sum_{j=1}^n dx_j\wedge dy_j$.
In \cite{minatlas} Y. B. Rudyak and F. Schlenk introduced the invariant $S_{B}(M,\omega)$ for a closed symplectic manifold $(M,\omega)$ of dimension $2n$ defined by:
\begin{equation*}
S_{B}(M,\omega):= \min  \lbrace  k \, \vert \, M= \mathcal{B}_{1}\cup\cdots\cup\mathcal{B}_{k}\rbrace,
\end{equation*}
where $\mathcal{B}_j$ is the image of a Darboux chart $\varphi(B^{2n}(r_j))\subset M$. This is the minimal number of symplectic charts needed to cover $(M,\omega)$. The problem of estimating this number is closely related to two other problems, namely computing the Gromov width $c_G(M,\omega)$ and the Lusternik-Schnirelmann category $\cat(M)$ of $M$.
While the latter can be often computed or estimated very well, computing the former is an open and delicate matter.
The Gromov width of a $2n$-dimensional symplectic manifold $(M, \omega)$, introduced in \cite{GROMOV85}, is defined as
\begin{equation*}\label{gromovwidth}
c_{G}(M,\omega)=\sup\big\lbrace\pi r^{2}\hspace{0.1cm}\big\vert\hspace{0.1cm}\exists\hspace{0.1cm}\varphi: \big(B^{2n}(r),\omega_{0}\big)\rightarrow(M,\omega)\big\rbrace
\end{equation*}
where $\varphi$ is a symplectic embedding.\\
By Darboux's theorem $c_G(M, \omega)$ is a positive number or $\infty$.
Computations and estimates of the Gromov width for various examples can be found in 
\cite{BIRAN97, BIRAN99, BIRAN01, castro, GROMOV85, JIANG00, GWgrass, LAMCSC11, GWHSSCT, LMZ02, LU06, LuDingQjao, MCDUFF91, MCDUFF94, SCHLENK05, GWcoadjoint}.\\
We adopt the following notation from \cite{GWHSSCT}.
\vskip 0.3cm
\noindent \label{pageNotation}
{\bf Notation:}  {\em From now on we shall use the shortening  HSSCT to denote a   Hermitian symmetric space of compact type.
Further, throughout the paper we shall denote by  $\omega_{FS}$   the canonical symplectic (\K) form on an irreducible HSSCT
normalized  so that $\omega_{FS} (B)\in\{-\pi, \pi\}$ when $B$ is a generator of $H_2(M, \Z)$, and by $A$ the generator for which $\omega_{FS} (A)= \pi$.}
\vskip 0.3cm


The following theorem and its two corollaries are the main results of this paper.

\begin{theor}\label{main00}
Let $(M,\omega_{FS})$ be a $2n$-dimensional HSSCT and let $f:M\hookrightarrow \CP^N$ be any holomorphic isometric immersion of $M$ in $\C P^N$ endowed with the Fubini--Study form $\omega$. Then
\begin{enumerate}[(i)]
\item If $\deg(f) \geq 2n$, then $S_{B}(M,\omega_{FS})=\deg(f)+1$
\item If $\deg(f) < 2n$, then $\max\lbrace n+1,\deg(f)+1\rbrace \leq S_{B}(M,\omega_{FS}) \leq 2n+1.$
\end{enumerate}
\end{theor}
As holomorphic isometric immersion $f:M\hookrightarrow\CP^N$ we can take, for example, the coherent states map described in Section \ref{The coherent states map}. In particular when $M$ is the complex Grassmannian one can take $f$ equal to the Pl\"ucker embedding. We recall the definition of degree of a holomorphic immersion in Section \ref{subs 1}, while in Section \ref{subs 2} we compute it for all irreducible HSSCT.

The proof of Theorem \ref{main00} is based on the results obtained by Y. B. Rudyak and F. Schlenk in \cite{minatlas} about minimal atlases for compact symplectic manifolds together with the explicit computation of the Gromov width given by the first author in collaboration with A. Loi and F. Zuddas in \cite{GWHSSCT} and the properties of the symplectic duality map introduced by  A. J.  Di Scala and A. Loi in \cite{DiScalaLoi08} which, in particular, give us a symplectic embedding of the noncompact dual $(\Omega, \omega_0)$ of  $(M, \omega_{FS})$ into  $(M, \omega_{FS})$.

Using the explicit computation of the volume of a classical domain $(\Omega, \omega_0)$ given by L. K. Hua in \cite{Hua}, we are able to prove the following corollary, which extends the computation of $S_{B}$ for the Grassmannians given in \cite{minatlas} to any classical irreducible HSSCT. Before stating the corollary, we recall that a classical irreducible HSSCT is one of the following quotients of compacts Lie groups:
\[
\I_{k,s}=SU(s)/S\left(U(k)\times U(s-k)\right),
\]
\[
\II_s=SO(2s)/U(s),
\]
\[
\III_s=Sp(s)/U(s),
\]
\[
\IV_s=SO(s+2)/SO(s)\times SO(2).
\]

\begin{cor}\label{main}
Let $(M,\omega_{FS})$ be a classical irreducible HSSCT of dimension $2n$. Then we have:
\begin{equation}\label{type1}
S_{B}( \I_{k,s})= \deg(f)+1,\quad \text { for }\quad (k=2 \text{ and } s \geq 7)\quad \text{or} \quad k \geq 3
\end{equation}
\[
S_B({\II_s})=\deg(f)+1, \quad \text { for }\quad  s \geq 6
\]
\[
S_B({\III_s})=\deg(f)+1, \quad \text { for }\quad  s \geq 5
\]
\[
n+1 \leq S_{B}(\IV_s) \leq 2n+1, \quad \text { for }\quad  s \geq 2.
\]
Otherwise, we have
\[
\max\lbrace n+1,\deg(f)+1\rbrace \leq S_{B}(M,\omega_{FS}) \leq 2n+1.
\]
\end{cor}
In the rank one case (i.e. $M=\C P^n$), we can set $f$ equal to the identity map, so that $\deg(f)=1$. On the other hand, \cite[Corollary 5.8]{minatlas} tells us that
\[
S_{B}(\C P^n,\omega_{FS})=n+1.
\]

The second corollary is a straightforward consequence of Theorem  \ref{main00}: 

\begin{cor}\label{main02}
Let $(M_1\times M_2,\omega_{FS})$ be a product of HSSCT of dimension $2n$. If $M_1\times M_2$ is different from $\C P^1 \times \C P^{n-1}$ and $\C P^2 \times \C P^2$, then
\begin{equation*}\label{mainformula}
S_{B}(M_1\times M_2,\omega_{FS})=\deg(f)+1,
\end{equation*}
where $f:M_1 \times M_2 \hookrightarrow \CP^N$ is any holomorphic isometric immersion.  Otherwise, we have $$\max\lbrace n+1,\deg(f)+1\rbrace \leq S_{B}(M,\omega_{FS}) \leq 2n+1.$$
\end{cor}

\noindent{\bf Acknowledgments.} The authors would like to thank Professor Andrea Loi for his help and various stimulating discussions and Professor Felix Schlenk for his interest in our work and his valuable comments. 

\subsection{The coherent states map.}\label{The coherent states map}
It is well know that an HSSCT $M$ is a simply connected \K--Einstein manifold with strictly positive scalar curvature. Therefore the integrality of $\frac {\omega_{FS}} {\pi}$ implies the existence of a polarizing holomorphic hermitian line bundle $(L,\, h)$ on $M$ such that $c_1(L)=[\frac {\omega_{FS}} {\pi}]$ and the Ricci curvature of $h$ satisfies $\Ric(h)=\frac {\omega_{FS}} {\pi}$ (where $\Ric(h)= -\frac{i} {2\pi} \, \de \deb \log \left(h\left(\sigma, \, \sigma\right)\right)$ in a local trivialization $\sigma: U \subset M \f L$). Consider the space $H^0(L)$ consisting of global holomorphic sections $s$ of $L$ which are bounded with respect to
\[
\langle s, s \rangle = \| s \| = \int_M h\left( s(x),\, s(x)\right) \frac{\omega^n}{n!}.
\] 

As $H^0(L) \neq \{0\}$, given an orthonormal basis $\{s_0, \dots, s_N\}\subset H^0(L)$ (with respect $\langle \cdot, \cdot \rangle$), it is well defined the \emph{coherent states map}, given by $f:M \f \C P^N$
\[
f(x)=\left[s_0(x): \dots : s_N(x)\right].
\] 
The Fubini--Study form $\omega$ of $\C P^N$ (normalized so that $\omega(B) \in \{-\pi,\pi\}$, when B is a generator of $H_2(\C P^N,\Z)$)
is given by $$\omega= \frac {i} {2} \, \de \deb \log \left(\sum_{j=0}^N \left| Z_j \right|^2\right),$$ 
it follows that
\[
f^*  {\omega}= \frac {i} {2}\, \de \deb \log \left(\sum_{j=0}^N \left| s_j(x) \right|^2\right)=\frac {i} {2} \,  \de \deb \log \left(\frac{\sum_{j=0}^N h\left( s_j(x),\, s_j(x) \right)}{h(\sigma(x), \, \sigma(x))}\right) 
\]
\[
=-\frac{i} {2} \, \de \deb \log \left(h(\sigma(x), \, \sigma(x))\right) + \frac{i} {2} \, \de \deb \log \left(\sum_{j=0}^N h\left( s_j(x),\, s_j(x) \right)\right) 
\]
\[
{\pi \,\Ric(h)}+  \frac {i} {2} \,  \de \deb \log \ep(x) = {\omega_{FS}} + \frac{i} {2} \,  \de \deb \log \ep(x),
\]
where $\ep:M \f \R$ is the so called \emph{$\ep$-function} defined by
\[
\ep(x)=\sum_{j=0}^N h\left( s_j(x),\, s_j(x) \right),
\]
one can prove that the $\ep$-function (see e.g. \cite[Theorem 4.3]{L00}) is invariant with respect the action of the group of holomorphic isometric transformation of $(M, \omega_{FS})$ which act transitively on $M$. Therefore the $\ep$-function is constant and we conclude that 
\[
f^* \omega=\omega_{FS}.
\]

\section{Proofs of Theorem \ref{main00}, Corollary \ref{main} and Corollary \ref{main02}}\label{sectionproofs}
Consider the following lower bound for $S_{B}(M,\omega)$ given by
\begin{equation*}\label{Gamma}
\Gamma(M,\omega):=\Bigg\lfloor\dfrac{Vol(M,\omega)n!}{c_{G}(M,\omega)^{n}}\Bigg\rfloor + 1,\hspace{0.3cm}
\end{equation*}
where $\lfloor x \rfloor$ denote the maximal integer smaller than or equal to $x$.
The following theorem summarizes the results about minimal atlases obtained in \cite{minatlas} that we need in the proof of Theorem \ref{main00}.
\begin{theorA}[Rudyak--Schlenk \cite{minatlas}]\label{teomin}
Let $(M, \omega)$ be a compact connected $2n$-dimensional symplectic manifold.
\begin{enumerate}
\item[i)]  If $\Gamma (M,\omega) \geq 2n+1$, then $S_{B}(M,\omega)=\Gamma(M,\omega)$.
\item[ii)] If $\Gamma (M,\omega) < 2n+1$ then $\max\lbrace n+1,\deg(f)+1\rbrace\leq S_{B}(M,\omega)\leq 2n+1$.
\end{enumerate}
\end{theorA}

\subsection{Proof of Theorem \ref{main00}}\label{subs 1}
We start recalling the definition of the degree of an holomorphic immersion $f\!:M\f\CP^N$. Suppose that $\dim(M)=2n<2N$, by Sard's Theorem there exists a point $q\notin f(M)$. Up to unitary transformation of $\CP^N$ we can suppose $q$ to be the point of coordinates $[1,0,\dots,0]$. Consider the projection $p_k\!:\CP^k\setminus \{q\}\f\CP^{k-1}$, $p_k([Z_0,\dots,Z_k])=[Z_1,\dots,Z_k]$ and define the map $F\!:M\f\CP^n$ by $F=\tilde p\circ f$, where $\tilde p=p_{n+1}\circ\cdots\circ p_N$. The degree $\deg(f)$ of $f$ is by definition the degree $\deg(F)$ of the map $F$, which is the integer number such that
\begin{equation}\label{deg}
F_*[M] =\deg (F)  [\CP^n] \in H_{2n}(\C P^n, \Z).
\end{equation}
What we need about $\deg(f)$ is summarized in the following Lemma:
\begin{lem}\label{lem4}{\rm(W. Wirtinger \cite{wirtinger}, M. Barros, A. Ros, \cite{bar-ros})} 
The degree $\deg (f)$ is a positive integer such that
\begin{equation}\label{voldeg}
\Vol(M)=\deg(f)\Vol(\C P^n),
\end{equation}
where $\deg (f)=1$ iff $M$ is totally geodesic and 
$\deg(f)=2$ iff  $f$ is congruent to the standard embedding of the quadric.
\end{lem}
The proof follows from Theorem A once one observes that the volume of any $n$-dimensional projective variety $X$, with holomorphic embedding $f:X\hookrightarrow\CP^N$, is given by
\begin{equation}\label{volume}
\Vol(X,\omega_{FS})=\deg(f)\Vol(\C P^{n}, \omega_{FS}),
\end{equation}
$\Vol(\C P^{n})=\frac{\pi^n}{n!}$ and that the Gromov width of any HSSCT (see \cite{GWHSSCT}) is given by $c_G(M,\omega_{FS})=\pi$.

\subsection{Proof of Corollary \ref{main}}\label{subs 2}

Consider $(\Omega, \omega_0)$, the noncompact dual of $(M,\omega_{FS})$.
In \cite[Theorem 1.1]{DiScalaLoi08} it is proved the existence of a global  symplectomorphism
\[\Phi: (\Omega, \omega_0)\rightarrow (M \setminus \Cut_p(M), \omega_{FS})\]
where $\Cut_p(M)$ is the cut locus of $(M, \omega_{FS})$ with respect to a fixed point $p\in M$ (see also \cite{diastexp}). Thus $\Vol(M,\omega_{FS})=\Vol(\Omega, \omega_0)$. On the other hand the explicit expression of the volume $\Vol(\Omega, \omega_0)$ can be found in L. K. Hua \cite{Hua} and by \eqref{volume} we are able to write the expression of $\deg(f)$ associated to any classical HSSCT, as follows.

\medskip

Let $\I_{k,s}$ be a HSSCT of type I, namely the Grassmannian of $k$-planes in $\C^s$. Notice that the dimension is $2n=2(s-k)k$ and that $\rank(\I_{k,s})=k$. We have that
\[
\deg(f_{k,s})=\frac{\Vol(\I_{k,s},\omega_{FS})}{\Vol(\C P^{(s-k)k}, \omega_{FS})}=
\]
\begin{equation}\label{deg grass}
=\frac{1!\,2! \dots (s-k-1)!\,1!\,2! \dots (k-1)! \left( (s-k)k \right)! }{1!\,2! \dots (s-1)!}.
\end{equation}
The case $\I_{k,s}$ was already done by Rudyak--Schlenk \cite{minatlas} and we obtain\eqref{type1} by \cite[Corollary 5.10]{minatlas}. Moreover they prove that
\[
S_B(\I_{2,4})\in\{5,6\}
\]
\[
S_B(\I_{2,5})\in\{7,8,9,10\}.
\]

Let $\II_s$ be an irreducible HSSCT of the second type. The complex dimension is given by $n_s=\frac{(s-1)s}{2}$. We have,
\[
\deg({f_{\operatorname{II_s}}})=\frac{s(s-1)}{2}!\, \frac{2! \, 4! \dots \left( 2s -4 \right)! }{\left(s-1\right)! \, s!  \dots \left(2s-3\right)!} .
\]
In order to apply Theorem \ref{main00} we need to study when
\[
\frac{\deg({f_{\operatorname{II_s}}})}{n_s}\geq 2
\]
One can see that the inequality is satisfied for $s=6$ and that $\frac{\deg({f_{\operatorname{II_s}}})}{n_s} < \frac{\deg({f_{\operatorname{II_{s+1}}}})}{n_{s+1}}$ for any $s\geq 6$.

Let $\III_s$ be an irreducible HSSCT of the third type. The complex dimension is given by $n_s=\frac{(s+1)s}{2}$. We have,
\[
\deg({f_{\operatorname{III_s}}})=\frac{s(s+1)}{2}!\,\frac{2! \, 4! \dots \left( 2s -2 \right)! }{s!\,\left(s+1\right)! \left(s+2\right)! \dots \left(2s-1\right)!  }.
\]
Arguing as before we see that $\frac{\deg({f_{\operatorname{III_s}}})}{n_s}\geq 2$ for any $s \geq 6$.

Let $\IV_s$ be an irreducible HSSCT of the fourth type (namely the complex quadric). Assume $s>3$ (if $s=1$ or $s=2$ we have respectively $\IV_1=\CP^1$ or $\IV_2=\CP^1\times \CP^1$).  By Lemma \ref{lem4}, $\deg(f)=2$. As $n=s \geq 3$, the result follows by $(ii)$ of Theorem \ref{main00}.

\subsection{Proof of Corollary \ref{main02}}
Let $\omega_{FS}^1$ and $\omega_{FS}^2$ be the Fubini-Study forms associated to $M_1$ and $M_2$. Since the associated volume form satisfies (with abuse of notation) $v_{\omega_{FS}}=v_{\omega_{FS}^1 } \wedge v_{\omega_{FS}^2}$, we have $\Vol(M_1\times M_2)=\Vol(M_1)\Vol (M_2)$. By (\ref{volume}) we get:
\begin{equation*}
\deg(f)=\frac{(n_1+n_2)!}{n_1! \, n_2! }\deg(f_1)\deg(f_2),
\end{equation*}
where $n_j$ is the complex dimension of $M_j$, $j=1,2$ and $f$, $f_1$ and $f_2$ are holomorphic isometric immersions of $M_1\times M_2$, $M_1$ and $M_2$.\ In order to apply $(i)$ of Theorem~\ref{main00}, we have to check when
\begin{equation}\label{degprod}
\deg(f_1)\,\deg(f_2)\, \frac{(n_1+n_2-1)!}{n_1! \, n_2! } \geq 2.
\end{equation}
First notice that when $\deg(f_1) \geq 2$ or $\deg(f_2) \geq 2$, since $\frac{(n_1+n_2-1)!}{n_1! \, n_2! } \geq 1$,  the inequality \eqref{degprod} is satisfied.

Assume now that $\deg(f_1)=\deg(f_2)=1$. By Lemma \ref{lem4}, $f_1,f_2$ are totally geodesic, this force $M_1$ and $M_2$ to have rank 1, that is $M_1=\C P^{n_1}$ and $M_2=\C P^{n_2}$. Moreover it is easy to see that  \eqref{degprod} is satisfied if and only if $n_1\geq 3$ and $n_2\geq 2$ or $n_1\geq 2$ and $n_2\geq 3$. The proof is complete.
\begin{remark}\rm
When $M = \C P^1 \times \C P^{n-1}, \C P^2 \times \C P^2$  we are not able to compute $S_{B}(M,\omega_{FS})$. Even for the simple case of  $\C P^1 \times \C P^1$ we know (private communication with F. Schlenk) that one can construct a covering
by  $4$ symplectic balls but we still do not know if this number can be reduced to $3$.\end{remark}


\begin{thebibliography}{99}\small


\bibitem{bar-ros} M. Barros, A. Ros, \emph{Spectral geometry of submanifolds}, Note di Matematica IV (1984), 1--56.


\bibitem{BIRAN97}
P. Biran,
\emph{Symplectic packing in dimension 4},
Geom. Funct. Anal. 7 (1997), 420-437.

\bibitem{BIRAN99}
 P. Biran,
 \emph{A stability property of symplectic packing},
 Invent. Math. 136 (1999) 123-155. 
 
 \bibitem{BIRAN01}
P. Biran, 
\emph{From symplectic packing to algebraic geometry and back}, 
European Congress of Mathematics, Vol. II (Barcelona, 2000)
Progr. Math. 202, Birkhauser, Basel (2001) 507-524.


\bibitem{castro}
A. C. Castro,
\emph{Upper bound for the Gromov width of coadjoint orbits of type A},  
arXiv:1301.0158v1

\bibitem{DiScalaLoi08} A. J.  Di Scala, A. Loi, \emph{Symplectic duality of symmetric spaces}, Advances in Mathematics 217 (2008), 2336-2352.

\bibitem{GROMOV85} M. Gromov, \emph{Pseudoholomorphic curves in symplectic manifolds}, 
Invent. Math. 82 (1985), no. 2, 307-347, Springer--Verlag (1986).

\bibitem{JIANG00}
M.-Y. Jiang,
\emph{Symplectic embeddings from $\R^{2n}$ into some manifolds},
Proc. Roy. Soc. Edinburgh Sect. A 130 (2000),  53-61.

\bibitem{Hua}
L. K. Hua,
\emph{Harmonic analysis of functions of several complex variables in the classical domains},
American Mathematical Society, Providence, R.I., 1963.



\bibitem{GWgrass}
Y. Karshon, S.   Tolman,
\emph{The Gromov width of complex Grassmannians},
Algebr. Geom. Topol. 5 (2005), 911-922.



\bibitem{LAMCSC11}
J. Latschev, D. McDuff and  F. Schlenk,
\emph{The Gromov width of $4$-dimensional tori},
Geom. Topol. 17 (2013) 2813-2853.


\bibitem{L00} C. Arezzo, Andrea Loi, 
\emph{Moment maps, scalar curvature and quantization of \K\ manifolds}, Comm. Math. Phys. 246(3) (2004), 543-559.


\bibitem{diastexp} A. Loi, R. Mossa,  {\em The diastatic exponential of a symmetric space}, 
 Math. Z. 268 (2011), no. 3-4, 1057-1068.
 
\bibitem{GWHSSCT} A. Loi, R. Mossa, F. Zuddas  {\em Symplectic capacities of Hermitian symmetric spaces of compact and non compact type}, to appear in J. Sympl. Geom.

\bibitem{LMZ02} A. Loi, R. Mossa, F. Zuddas,
\emph{Some remarks on the Gromov width of homogeneous Hodge manifolds}, Int. J. Geom. Methods Mod. Phys. 11 (2014), no. 9.

\bibitem{LU06}
G. Lu,
\emph{Gromov-Witten invariants and pseudo symplectic capacities}
 Israel J. Math. 156 (2006), 1-63.

\bibitem{LuDingQjao}
G. Lu, H. Ding, Q.  Zhang,
\emph{Symplectic capacities of classical domains},  
Int. Math. Forum 2 (2007), no. 25-28, 1311-1317. 


\bibitem{MCDUFF91}
D. McDuff, \emph{Blowups and symplectic embeddings in dimension 4},
Topology 30 (1991), 409-421.

\bibitem{MCDUFF94}
D. McDuff and L. Polterovich,
\emph{Symplectic packings and algebraic geometry},
Invent. math. 115 (1994), 405-429.


\bibitem{minatlas}
Y. B. Rudyak, F. Schlenk,
\emph{Minimal atlases of closed symplectic manifolds},  
Commun. Contemp. Math. 9 (2007), no. 6, 811-855. 

\bibitem{SCHLENK05}
F. Schlenk,
\emph{Embedding problems in symplectic geometry},
de Gruyter Expositions in Mathematics 40. Walter de Gruyter Verlag, Berlin, 2005.

 


\bibitem{ta} M. Takeuchi,
{\em Homogeneous \K\ Submanifolds in Complex Projective Space}, Japan
J. Math. vol. 4 (1978), 171-219.


\bibitem{wirtinger} W. Wirtinger, \emph{Eine Determinantenidentit\"at und ihre Anwendung auf analytishe Gebilde in Euclidischer und Hermitischer Massbestimmung} Monatsh. Math. Phys. 44 (1936), 343--365.





\bibitem{GWcoadjoint}
M. Zoghi, 
\emph{The Gromov Width of Coadjoint Orbits of Compact Lie Groups}, 
Thesis (Ph.D.) University of Toronto (Canada)  2010.





















\end{thebibliography}
\end{document}